\documentclass{article}

\begin{document}
\large
\title{A Question about Differential Ideals}         

\author{Eloise Hamann\\Department of Mathematics,\\San Jose State University,San Jose, Ca 95192-0103\\hamann@math.sjsu.edu}
 
\maketitle
\textbf{Abstract:} This paper is an extract of a paper that will soon appear in Communications in Algebra.  The paper investigates the converse to a theorem which states that a differential domain D finitely generated over a field F of characteristic 0 with no nonzero proper differential ideals will have constants of differentiation C only in the field in question.  The converse is known to be false but the question of whether the differential domain D can be finitely extended within its quotient field E to a domain with no nonzero proper differential ideals was raised in [M]. The paper about to appear contains two examples one of which contains a fatal error.  This extract contains one example to show that the answer to the question posed is no.  It is shown to be true under hypotheses which include restrictions on F and limiting the number of generators of D over C to two or  the transcendence degree of D over C to one.  \\

\section{Introduction}       
Consider the following theorem about differential algebras.\\

\textbf {Theorem:}  Let R be a differential integral domain, finitely generated over the differential field F of characteristic 0.  Let E denote the quotient field of R, and let C denote the field of constants of F.  Assume that R contains no proper differential ideals and that the field of constants C of F is algebraically closed of characteristic 0.  Then the field of constants of E coincides with C.  [M]  \\

Andy Magid shows in [M] that the converse is false, namely an example is given of a differential domain R, finitely generated over the field of constants C of F where the field of constants of E coincides with C (a situation described by saying that R has no new constants), but R has proper nonzero differential ideals. The R in the example does have a finitely generated extension contained in its quotient field which has no nonzero proper differential ideals so the question is raised  whether it is always possible to find a finitely generated F-subalgebra T of E containing R such that T has no proper differential ideals when R has "no new constants". If this were the case one has the sense that the property of having no new constants would be governed by not having too many nonzero proper differential ideals. \\ 

  All fields in this paper are assumed to have characteristic 0. We assume basic familiarity with differential ideals (ideals closed under differentiation), differential homomorphisms (those that commute with differentiation) whose kernels are differential ideals, and the correspondence theorem for differential rings which are homomorphic images of a differential ring.  We include a lemma that is described as folklore.\\  

\textbf{Lemma 0.1}  In a Noetherian ring R containing the rationals, a minimal prime P over a differential ideal I is differential.\\

\textbf{Pf:}  In any ring containing the rationals, the radical of a differential ideal is known to be differential.  See 4.2 of [L].  Let J be the radical of a differential ideal I.  If P is a minimal prime over I, it is also a minimal prime over J and since R is Noetherian, J = Ann\{x\} for some x $\notin$ J.  If y $\in$ P, yx $\in$ J so xDy + yDx $\in$ J.  Then x$^{2}$Dy $\in$ J so (xDy)$^{2}$ $\in$ J, but J is radical so xDy $\in$ J and Dy $\in$ P so P is differential.

\section{Counterexample}

\textbf{Example}  Let F = C(\{$\alpha_{i}$\}), i $\in$ Z$^{+}$ where $\alpha_{i}$ are independent indeterminates.  Define a derivation D on F by Dc = 0 if c $\in$ C and D($\alpha_{i}$) = $\alpha_{i}^{3} - 2\alpha_{i}^{2} + 2\alpha_{i}$ and extend to F by the rules of derivation.  Let R = F[X], E the quotient field of R and extend D to E by defining D(X) = X$^{3} - 2X^{2} +2X$.  In the following we will call elements of C, C-constants and elements w such that Dw = 0, differential constants. We will refer to the $\alpha_{i}$ as variables. \\

\textbf{Claims:}  

\textbf{1)} The differential constants of F under D = the C-constants, C.

\textbf{2)} The differential constants of E also = C.

\textbf{3)} $($X $-$ $\alpha_{i}$$)$ is a differential ideal of R.  

\textbf{4)} Every finitely generated extension of R which is contained in E has a proper differential ideal.       

\textbf{Pf of 1)} If w $\in$ C[\{$\alpha_{i}$ \}] is of positive degree n in some $\alpha_{i}$, then it is easy to check that Dw has degree n + 2 in $\alpha_{i}$.  Thus, Dw $\neq$ 0 when w is not a C-constant in the polynomial ring,C[\{$\alpha_{i}$ \}].  Let w $\in$ F and write w = f/g with f and g relatively prime polynomials in C[\{$\alpha_{i}$ \}]. If w is a differential constant, then Dw = 0 and gDf $-$ fDG = 0.  Since f and g are relatively prime, we must have Df = zf and Dg = zg for some z in C[\{$\alpha_{i}$ \}].  Thus, it suffices to show that given z and a particular nonzero solution f to DY = zY, any other solution is a C-constant multiple of f.   Since applying D increases the degree of each $\alpha_{i}$ which occurs by two, z is of degree two in each $\alpha_{i}$ which occurs in f.\\
Renumber so that f involves $\alpha_{i}$, i = 1,...,n. We proceed by induction on n so let f = $\sum_{i=1}^N a_{i}\alpha_{1}^{i}$ with $a_{i}\in$ C.  Df =$\sum_{i=1}^N ia_{i}\alpha_{1}^{i+2} - 2\sum_{i=1}^N ia_{i}\alpha_{1}^{i+1} + 2\sum_{i+1}^N ia_{i}\alpha_{1}^{i}$ and fz = $\sum_{i=1}^N a_{i}\alpha_{1}^{i}$(b$_{2}\alpha_{1}^{2}$ + b$_{1}\alpha_{1}$ + b$_{0}$).  Equating coefficients we obtain:  

N$a_{N}$ = b$_{2}a_{N}$ so b$_{2}$ = N.

(N $-$ 1)$a_{N-1} - 2Na_{N} = Na_{N-1} + a_{N}$b$_{1}$ so $a_{N-1} = -$(b$_{1} +2N)a_{N}$.

For all i $\le$  N $-$ 2, $ia_{i} -2(i+1)a_{i+1} +2(i+2)a_{i+2} = Na_{i}$ + b$_{1}a_{i+1}$ + b$_{0}a_{i+2}$ so $(N - i)a_{i} = (2i + 4 -$ b$_{0})a_{i+2} -$ (b$_{1} + 2(i+1))a_{i+1}$.  

Since $a_{N} \neq 0$ and we are only arguing uniqueness up to a nonzero multiple we may assume $a_{N}$ = 1. The relations imply that b$_{2}$ = N, $a_{N-1}$ is determined by b$_{1}$ and recursively downward all remaining $a_{i}$ are determined by b$_{0}$ and b$_{1}$ and previously determined $a_{i}$.  Thus, the claim follows for n = 1.  

Now suppose true for fewer than n variables with n $\ge$ 2 and that f is a nonzero solution to DY = zY involving n variables so f = $\sum_{i=1}^{N}f_{i}\alpha_{n}^{i}$ with $f_{i}\in$ C[$\alpha_{1},\ldots,\alpha_{n-1}$].  We may assume that f$_{0}\neq$ 0 since if f = $\alpha_{n}^{j}$g where $\alpha_{n}$ does not divide g, then Df = zf implies that Dg = (z $- j\alpha_{n}^{2} + 2j\alpha_{n} - 2j)g$ so if g is determined up to a nonzero constant multiple then clearly so is f. We have:

\textbf{(1)} Df = $\sum_{i=1}^N Df_{i}\alpha_{n}^{i}$ + $\sum_{i=1}^N if_{i}\alpha_{n}^{i+2} - 2\sum_{i=1}^N if_{i}\alpha_{n}^{i+1} + 2\sum_{i+1}^N if_{i}\alpha_{n}^{i}$ and fz = $\sum_{i=1}^N f_{i}\alpha_{n}^{i}$(b$_{2}\alpha_{n}^{2}$ + b$_{1}\alpha_{n}$ + b$_{0}$) with b$_{i}\in$ C[$\alpha_{1},\ldots,\alpha_{n-1}$].        

As in the n = 1 case, comparing leading coefficients gives b$_{2}$ = N.  We also  have Df$_{0}$ = f$_{0}$b$_{0}$.  By induction, f$_{0}$ is determined up to a C-constant multiple by b$_{0}$.  In addition:

Df$_{1}$ + 2f$_{1}$ = b$_{0}f_{1}$ + b$_{1}f_{0}$ so Df$_{1}$ = (b$_{0} - 2$)f$_{1}$ + b$_{1}f_{0}$.  

For i $\ge$ 2, Df$_{i}$ + (i $-$ 2)f$_{i-2}$  $-$ 2(i $-$ 1)f$_{i-1}$ + 2if$_{i}$ = b$_{0}f_{i} +$ b$_{1}f_{i-1} + Nf_{i-2}$ so Df$_{i}$ = (b$_{0} - 2i$)f$_{i}$ + (b$_{1}$ + 2(i $-$ 1))f$_{i-1}$ + (N $-$ i + 2)f$_{i-2}$.  

In all cases Df$_{i}$ =(b$_{0} - 2i$)f$_{i}$ + P$_{i}$ where P$_{i}$ involves earlier f$_{i}$.  Assume the following *Lemma for n variables.

\textbf{*Lemma:}  Consider F$_{n}$ = C($\alpha_{1},\ldots,\alpha_{n}$).  DY = cY where c is a nonzero C-constant has only the trivial solution Y = 0 in F$_{n}$.  

Then suppose we have another solution, g,  to DY = zY. \textbf{Case 1.}  g$_{0}$ = 0.  Write g = $\alpha_{n}^{j}$h where j $\ge$ 1 and $\alpha_{n}$ does not divide h.  Then we obtain Dh = (z $-j\alpha_{n}^{2} +2j\alpha_{n} -$ 2j)h. The same argument that gives Df$_{0}$ = b$_{0}f_{0}$ gives Dh$_{0}$ = (b$_{0} -$ 2j)h$_{0}$.  We have h$_{0}\neq$ 0 and D(f$_{0}/h_{0}$) = 2j(f$_{0}/h_{0}$) contradicting the *Lemma for n $-$ 1 variables.  \textbf{Case 2.} g$_{0} \neq$ 0.  Now Dg$_{0}$ = b$_{0}g_{0}$ so by induction g$_{0}$ = cf$_{0}$ with c $\neq$ 0.  Since (1/c)g is a solution to DY = zY if g is a solution we may assume that g$_{0}$ = f$_{0}$.  Now g$_{i}$  must satisfy exactly the same recursion relations satisfied by f$_{i}$.  In particular, Dg$_{1}$ = (b$_{0} -$2)g$_{1}$ + b$_{1}g_{0}$.  Since g$_{0}$ = f$_{0}$, f$_{1}$ and g$_{1}$ are both solutions of DY = (b$_{0} -$ 2)Y   + b$_{1}$f$_{0}$.  Therefore g$_{1} - f_{1}$ is a solution  of DY = (b$_{0} -$ 2)Y and (g$_{1} - f_{1}$)/f$_{0}$ is a solution of DY = $-$ 2Y as in the proof of Case 1.  The *Lemma for n $-$ 1 variables then implies that g$_{1} - f_{1}$ = 0.  Now assume that f$_{i}$ = g$_{i}$ for 0 $\le$ k where k $\ge$ 1.  Since the starting two coefficients of f and g agree, the common recursion relation gives that both g$_{k+1}$ and f$_{k+1}$ are solutions of DY = (b$_{0} -$ 2(k+1))Y + P$_{k+1}$ so that (g$_{k+1} - f_{k+1}$)/f$_{0}$ is a solution of DY = $-$ 2(k+1)Y so by the *Lemma for n $-$ 1 variables, g$_{k+1}$ = f$_{k+1}$ and by induction g = f.  

It remains to prove the *Lemma. Again we use induction on the number of variables. It should first be noted that the argument above shows that if the *Lemma holds for n $-$ 1 variables or fewer that Df = zf with f in n variables implies that z determines a nonzero solution f up to a C-constant multiple of f if such an f exists.  Suppose that w($\alpha_{1}$) = f($\alpha_{1}$)/g($\alpha_{1}$) is a solution of DY = cY in the case of one variable where c is a nonzero C-constant.  Then w$'$($\alpha_{1}$)($\alpha_{1}^{3} -$2$\alpha_{1}^{2}$ + 2$\alpha_{1}$) = cw($\alpha_{1}$).  c and w are $\neq$ 0, and this is an ordinary calculus problem whose solution is w = A$\alpha_{1}^{c/2}$[($\alpha_{1} - 1)^{2}$ + 1]$^{-c/4}$e$^{(c/2)tan^{-1}(\alpha_{1} - 1)}$ with A $\neq$ {0} which is not a rational function of $\alpha_{1}$ unless c = 0.  If w were rational then in the case that c is rational every power of w is also rational and one could conclude that e$^{ktan^{-1}(\alpha_{1} - 1)}$ would be rational for some k $\neq$ 0 which is impossible.  If c is irrational, then one could conclude that e$^{(c/2)tan^{-1}(\alpha_{1} -1)})$ is either 0 or undefined at $\alpha_{1}$ = 0 again a contradiction. Thus the *Lemma holds for n = 1.

Now suppose the *Lemma is true for fewer than n variables and that D(f/g) = c(f/g) with c $\neq$ 0. We assume neither f nor g are 0 and are relatively prime. Write f = $\sum_{k=0}^N f_{k}\alpha_{n}^{k}$ and g = $\sum_{k=0}^N g_{k}\alpha_{n}^{k}$.     We have D(g/f) = $-$c(g/f) so if f/g is a counterexample so is g/f. WLOG assume that f$_{0}\neq$ 0.   [gDf $-$ fDg]/g$^{2}$ = c(f/g) so gDf $-$ fDG = cfg and f must divide Df and g must divide Dg.  Let Df = zf, then zgf $-$ fDg = cfg and Dg = (z $-$ c)g. It is easily checked that D commutes with substituting $\alpha_{n}$ = 0 and in fact that D commutes with substituting $\alpha_{n}$ = 1 $\pm$ $i$ as all are roots of D$\alpha_{n}$.  Thus, if g$_{0} \neq$ 0, D(f$_{0}$/g$_{0}$) = c(f$_{0}$/g$_{0}$) and we contradict the induction hypothesis.  Thus, write g = $\sum_{k=j}^N g_{k}\alpha_{n}^{k}$ where j $>$ 0 and g$_{j}\neq$ 0.  Writing z as N$\alpha_{n}^{2}$ + b$_{1}\alpha_{n}$ + b$_{0}$, Df = zf as in (1) above implies that Df$_{0}$ = b$_{0}$f$_{0}$. A similar argument implies that Dg$_{j}$ + 2jg$_{j}$ = (b$_{0}$ $-$ c)g$_{j}$ so Dg$_{j}$ = g$_{j}$(b$_{0}$ $-$c $-$ 2j).  If c + 2j $\neq$ 0, then D(f$_{0}$/g$_{j})$ = (c + 2j)(f$_{0}$/g$_{j}$) and the induction hypothesis would imply that f$_{0}$ = 0 providing a contradiction.  Thus, we must have c = $-$2j and Dg$_{j}$ = g$_{j}$b$_{0}$.   Since f$_{0}$ and g$_{j}$ involve n $-$ 1 variables by the remark in the opening of the previous paragraph, we must have g$_{j}$ = kf$_{0}$ for some k $\neq$ 0.  Since (1/k)g must also satisfy DY = (z + 2j)Y if g does, we may assume that g$_{j}$ = f$_{0}$ and D(f/g) = $-$2j(f/g).   

Now we repeat the above argument expanding f by $\alpha_{n} -$ (1 + $i$), writing f = $\sum_{k=1}^N$ f$_{k}$*($\alpha_{n} - i -$ 1)$^{k}$ and  obtain a different value of c than 2j.  Recall that D$\alpha_{n}$ = $\alpha_{n}^3 - 2\alpha_{n}^2 + 2\alpha_{n}$ =$\alpha_{n}$($\alpha_{n} - i -$ 1)($\alpha_{n} + i -$ 1). Let $\alpha_{n}$* = $\alpha_{n} - i -$ 1.  Both f and g cannot be divisible by $\alpha_{n}$*. We now assume that we only know that D{f/g} = $\pm$ 2j i.e. we don't know which of f or g is not divisible by $\alpha_{n}$*.  Since D(f/g) = c(f/g) implies that D(g/f) = $-$c(g/f) in arguing that expansion about $\alpha_{n}$* contradicts the c value being $\pm$ 2j there is no loss of generality in assuming that the numerator polynomial is f and that f$_{0} \neq$ 0. In terms of $\alpha_{n}$*, z = N$\alpha_{n}$*$^{2} + b_{1}$*$\alpha_{n} + b_{0}$* where b$_{1}$* = 2N +2N$i$ + b$_{1}$ and b$_{0}$* = 2N$i$ + b$_{1}$(1 + $i$) + b$_{0}$ although the actual values are not relevant. Let g = $\sum_{k=r}^N g_{k}$*$\alpha_{n}$*$^{k}$ with g$_{r}$* $\neq$ 0.  D$\alpha_{n}$* = D$\alpha_{n}$ = $\alpha_{n}^{3} - 2\alpha_{n}^{2} + 2\alpha_{n}$ = $\alpha_{n}$*$^{3} +$ (3$i$ + 1)$\alpha_{n}$*$^{2}$ + (2$i$ $-$ 2)$\alpha_{n}$*.  If g$_{0} \neq$ 0, then setting $\alpha_{n}$* = 0 which commutes with D yields D(f$_{0}$*/g$_{0}$*) = $\pm$\,2j(f$_{0}$*/g$_{0}$*) contradicting the induction assumption.  Thus, g$_{0}$ = 0  and r $>$ 0.

By a nearly familiar argument, we can derive Df$_{0}$* = b$_{0}$*f$_{0}$* from expressing Df = zf in terms of $\alpha_{n}$* below:

$\sum_{k=0}^N Df_{k}$*$\alpha_{n}$*$^{k} + \sum_{k=0}^N kf_{k}$*$\alpha_{n}$*$^{k+2}$ + (3$i$ + 1)$\sum_{k=0}^N kf_{k}$*$\alpha_{n}$*$^{k+1}$ + (2$i$ $-$ 2)$\sum_{k=0}^N kf_{k}\alpha_{n}$*$^{k}$ = $\sum_{k=0}^N f_{k}$*$\alpha_{n}$*$^{k}$(N$\alpha_{n}$*$^{2}$ + b$_{1}$*$\alpha_{n}$* + b$_{0}$*).

Since D(f/g) = $\pm$2j(f/g),  Dg = (z $\mp$ 2j)g, and  we have: 

$\sum_{k=r}^N Dg_{k}$*$\alpha_{n}$*$^{k} + \sum_{k=r}^N kg_{k}$*$\alpha_{n}$*$^{k+2}$ + (3$i$ + 1)$\sum_{k=r}^N kg_{k}$*$\alpha_{n}$*$^{k+1}$ + (2$i$ $-$ 2)$\sum_{k=r} kg_{k}$*$\alpha_{n}$*$^{k}$ = $\sum_{k=r}^N g_{k}$*$\alpha_{n}$*$^{k}$(N$\alpha_{n}$*$^{2}$ + b$_{1}$*$\alpha_{n}$* + b$_{0}$* $\mp$ 2j).

We obtain Dg$_{r}$* = (b$_{0}$*$\mp$2j  $-$  (2$i$ $-$ 2)r)*g$_{r}$*.  Thus, D(f$_{0}$*/g$_{r}$*) = [$\mp$2j $-$ (2$i$ $-$ 2)r](f$_{0}$*/g$_{r}$*).  This contradicts the induction assumption unless $\mp$2j $-$ (2$i$ $-$ 2)r = 0 which is impossible since j and r are not 0. Thus, (1) holds.

\textbf{Pf of (2):}  This is clear since E is differentially isomorphic to F.  

\textbf{Pf of (3):}  This follows from D(X $-$ $\alpha_{i}$) = X$^{3}$ $-$ 2X$^{2}$ + 2X $-$ ($\alpha_{i}^{3} - 2\alpha_{i}^{2} + 2\alpha_{i}$) = (X $-$ $\alpha_{i}$)(X$^{2} + \alpha_{i}$X +  $\alpha_{i}^{2} -$ 2(X + $\alpha_{i}$) + 2).  

\textbf{Pf of (4):}  It is clear that the generators of any finitely generated extension of R in E can involve only finitely many $\alpha_{i}$. Thus, an ideal generated by X $- \alpha_{j}$ must remain proper for some j and clearly remains a differential ideal. \\ 

\textbf{Remark:} The example above has a defining polynomial, X$^{3} - 2X^{2} +2X$.  The author investigated simpler polynomials which did not work.  However, any polynomial p(X) in Z[X] where DY/Y = DX/(p(X)) does not have a rational solution Y = r(X) and p(X) has  a pair of roots linearly independent over Q should suffice. Thus, the counterexample is one member of an infinite family of counterexamples.

Recall the general setting of the question where C $\subseteq$ F $\subseteq$R.  In the example above, R was generated over F simply by a single element.  However, F was not finitely generated over C. Thus, the example fails to be geometric in some sense.   It is natural to ask if requiring F to be finitely generated over C can produce a positive answer. The paper to appear has an example purporting to answer the question in the negative. The primary reason for posting this extract is to alert readers that this question remains open.

\section{Positive Results}

Since we will be working with polynomial rings, the term differential constant will again be occasionally employed to avoid confusion.  \\

\textbf {Def:}  Let R be a differential ring with derivation D.  Call a prime differential ideal P an integral ideal if gDf $-$ fDG $\in$ P implies that f $-$ cg or g $-$ cf $\in$ P for some differential constant  c.  \\

\textbf {Observation:}  If P is an integral ideal in a differential ring then Df $\in$ P implies  that f $-$ c $\in$ P for some differential constant c under D.\\  

\textbf {Pf:} Let g = 1, then 1Df $-$ fD1 $\in$ P so f $-$ c$\cdot$1 $\in$ P or 1 $-$ fc $\in$ P. In the latter case c $\neq$ 0 since P is proper.  Then f $-$ 1/c $\in$ P.\\  

\textbf {Lemma 1.}  Assume n is a fixed integer and suppose the following for all  polynomial rings R = C[X$_{1}$,...,X$_{n}$] equipped with a derivation D,  where Dc = 0 if c $\in$ C and C is any field of characteristic 0.  Given any integral ideal P of R  either P is a maximal differential ideal or the intersection of all prime differential ideals properly containing P also properly contains P.  Then we can conclude that all differential domains D which can be generated by n elements over the algebraically closed field C where C is the field of differential constants of the quotient field, E,  of D have finitely generated extensions within E which have no nonzero proper differential ideals.  \\

\textbf{Pf:} Let  C = the  field of differential constants of E, the quotient field of D = C[x$_{1},\ldots$,x$_{n}$].  Let D be the obvious homomorphic image of R = C[X$_{1},\ldots$,X$_{n}$].  Define DX$_{i}$ to be a polynomial of least degree such that DX$_{i}$ maps to Dx$_{i}$. It is easy to see that R is a differential homomorphic image of T.  Let P = Ker $\theta$ where $\theta$ is the map from T onto R.  It is clear that P is a differential prime ideal.  If gDf $-$ fDg $\in$ P in T, then $\theta$(g)D$\theta$(f) $-$ $\theta$(f)D$\theta$(g) = 0 so since the differential constants of R = C, we must have either both $\theta$(f)and $\theta$(g) are 0 or $\theta$(f)/$\theta$(g) = c, or $\theta$(g)/$\theta$(f) = c for some c $\in$ C.  Thus, f $-$ cg $\in$ P or  g $-$ cf $\in$ P for some c so P is an integral ideal.  By hypothesis, either P is not properly contained in any prime differential ideal or the intersection of such, properly contains P.  In the former case, P is a maximal differential ideal so that R has no nonzero proper differential ideals.  In the latter case, there is an element, w,  in the intersection of all prime differential ideals containing P which is not in P.  Thus, $\theta$(w) is not 0, but is in all nonzero prime differential ideals of R by a correspondence theorem for differential ideals.  Clearly R[1/$\theta$(w)] is a finitely generated extension of R within E which has no nonzero prime differential ideals.\\                 

\textbf{Lemma 2:}  Let R = C[X$_{1}$,\ldots,X$_{n}$] be a polynomial ring equipped with a derviation D, with C the field of constants of the quotient field E of R.  Suppose C is finite algebraic over Q. Then  R has only finitely many height one differential prime ideals.\\  

\textbf{Pf:}  Let P be a height one differential prime ideal.  Since R is a UFD, P is principal and if f generates P we must have Df = wf for some w $\in$ R.  Let V = \{w $\in$ R $\mid \exists$ f so Df = wf for some f\}.  Let d = maximum total degree of the DX$_{i}$.  Let B be a finite basis of C over Q. 

Let U be the Q vector space generated by all products of elements of B with a monomial in \{X$_{i}$\} where each variable occurs with degree at most d.  Let M be this finite set of monomials.  By construction U is a finite dimensional vector space over Q.  

We claim that V $\subseteq$ U.  Let w $\in$ V and f be in C[X$_{1}$,\ldots,X$_{n}$]  such that Df = wf.   Since Df = $\sum_{i=1}^{n}\partial$f/$\partial$X$_{i}$DX$_{i}$, it is clear that the total degree of f in the \{X$_{i}$\} is at most degree f $-$ 1 + d so that w is in the C vector space generated by M and since C is generated over Q by B, w $\in$ V.    

Since U is a finite dimensional vector space over Q, there exist at most a finite number of w in V which are linearly independent over Q.  Let A = \{w$_{1}$,\ldots,w$_{N}$) be a maximal linearly independent subset of V that has the property that the corresponding f$_{i}$ are irreducible.    

Now we show that \{P$_{i}$ = (f$_{i}$)\} is the entire set of height one differential primes of R.  Suppose P is such an ideal and that P = (g).  Thus, g is irreducible and Dg = wg so w $\in$ V.  Then we have since A is a maximal linearly independent set over Q corresponding to irreducibles that mw = $\sum_{i=1}^{N}$n$_{i}$w$_{i}$ for integers m and \{n$_{i}$\}.  Then Dg$^{m}$ = mwg$^{m}$.  We also have that if z = $\prod_{i=1}^{N}$f$_{i}^{n_{i}}$ that Dz = $\sum_{i=1}^{N}$n$_{i}$w$_{i}$ z = mwz.  However, we now have that D(g$^{m}$/z) = (zDg$^{m}$ $-$ g$^{m}$Dz)/z$^{2}$ = 0 so g$^{m}$/z $\in$ C. g$^{m}$/z = c implies that g$^{m}$ = zc.  z is a product of powers of the irreducbile {f$_{i}$}.  Some of the powers may be negative but in any case by unique factorization, g must be an associate of one of the f$_{i}$.  Thus, P =(g) = (f$_{i}$) and the Lemma holds.\\

The following theorem gives a positive answer when two elements generate D over a field C finite algebraic over Q. Thus, it includes the case of a vector field on the plane, C$^2$.

\textbf{Theorem 3.} Let D = C[x$_{1}$,x$_{2}$] be a domain equipped with a derivation D,  where the differential field of constants of E, the quotient field of D, is C. Suppose C is finite algebraic over Q.   Then D has a finitely generated extension within E with no nonzero proper differential ideals.\\  

\textbf {Pf:}  By Lemma 1. it suffices to show that the integral ideals of any  polynomial ring R = C[X$_{1}$, X$_{2}$] with a derviavtion D are either maximal or are not intersections of the prime differential ideals properly containing them.  Consider first the case where (0) is an integral ideal.  By Lemma 2. there exists g $\neq$ 0  $\in \cap$ P$_{i}$, the intersection of the finite number of differential height one primes. Let Q* be the algebraic closure of Q and extend D to Q*[X$_{1}$, X$_{2}$].  Dc = 0 if c $\in$ Q*.  Each maximal ideal M of D has height 2 and its extension to D* = Q*[X$_{1}$, X$_{2}$] also has height 2 so any maximal ideal of D* which contains M is minimal.  It is easy to check that M extends to a differential ideal of D*.  By 0.1, a minimal prime over the extension of M is differential.  Since maximal ideals of D* are generated by \{X$_{1}-$ a, X$_{2} -$ b\} for some constants a and b, a maximal differential ideal of D* must contain both DX$_{i}$ neither of which is 0 since (0) is an integral ideal of D. But each maximal ideal of D* must intersect D in M. Thus, DX$_{1}$ $\in$ M for every differential max ideal M of D so gDX$_{1}$ is in every proper nonzero differential prime ideal.  Thus, 0 is not the intersection of differential prime ideals properly containing 0.  Now suppose that P is a height one integral ideal. Say P = (h).   Let DX$_{1}$ = yz and DX$_{2}$ = wz where w and y are relatively prime.  Since P is an integral ideal z $\notin$ P since then DX$_{i}$ would be in P so X$_{i}$ + c$_{i}$ would be in P for some c$_{i}$ for i = 1,2 so that P would be maximal. Any differential maximal ideal  either contains both y and w or it contains z by arguing in D*.  Since h must be irreducible and y and w are relatively prime, one of y or w $\notin$ P, say y$\notin$ P.  Then we have yz is not in P but yz $\in \cap$ M such that M is a maximal differential ideal containing P.  The claim follows.\\  

The following theorem is evidently well known as it can be interpreted as saying that given a non-zero vector field on an algebraic curve there will be a Zariski open set on which the vector field doesn't vanish.  A proof is given in the spirit of this paper.\\  

\textbf{Theorem 4.} Let C be an algebraically closed field of characteristic 0.   If D = C[{x$_{1}$,\ldots,x$_{n}$] is a domain equipped with a derivation D, of transcendence degree one over C or equivalently of Krull dimension one, where C is the field of differential constants of E, the quotient field of D,  then D has a finitely generated extension within E which has no nonzero proper differential ideals.\\ 
 
\textbf{Pf:}  Define D on  R = C[{X$_{1}$,\ldots,X$_{n}$] by defining DX$_{i}$ to be a particular lift of Dx$_{i}$.  If $\theta$ is the obvious map from R to D, ker $\theta$ is an integral ideal of height n $-$ 1 since R has Krull dimension = the transcendence degree of R over C = 1.  As in the proof of Lemma 1, it suffices to prove that the intersection of differential prime ideals properly containing ker $\theta$ is not ker $\theta$.  Since ker $\theta$ has height n $-$ 1, any such prime differential ideal must be maximal.  Because ker $\theta$ is an integral ideal which is nonmaximal, some DX$_{j}$ $\notin$ ker $\theta$. Since C is algebraically closed, any maximal ideal is generated by n elements of the form X$_{i}$ $-$ c$_{i}$.  Then all DX$_{i}$ must be contained in a differential maximal ideal.  Thus, DX$_{j}$ is in the intersection of differential primes containing ker $\theta$, but not in ker $\theta$.   
\begin{center}
\textbf{Bibliography}
\end{center}

[K1] I. Kaplansky,An Introduction to Differential Algebra, Paris. Hermann 1957.

[K2] I. Kaplansky, Commutative Rings, Allyn and Bacon, Boston. 1970
  
[L] T. Little, Algebraic Theory of Differential Equations, Thesis San Jose State University, 1999

[M]A. Magid, Lectures on Differential Galois Theory, University Lecture Series V.7, Providence, AMS, 1991

\end{document}